\documentclass[a4paper,12pt]{article}
\usepackage{localeng}
\usepackage[dvipsnames]{xcolor}
\usepackage{skak}
\usepackage{url}
\usepackage[margin=25mm]{geometry}
\newcounter{probcounter}
\newenvironment{pr}%
   {\medskip\par\refstepcounter{probcounter}\color{Sepia}\textbf{Example \arabic{probcounter}.}  }%
   {\par\medskip}%

\begin{document}
\title{Constructive mathematics and teaching}
\author{Alexander Shen\thanks{LIRMM, University of Montpellier, CNRS, Montpellier, France, \texttt{alexander.shen@lirmm.fr}, \texttt{sasha.shen@gmail.com}. This note was prepared for the invited talk at the \textsc{Constructive Mathematics: Foundations and Practice, CM:FP 2023} conference in memory of Errett Bishop, June 2023. I~am grateful to the organizers of this conference for their invitation and support.}}
\maketitle

\begin{abstract}
Constructivists (and intuitionists in general) asked what kind of mental construction is needed to convince ourselves (and others) that some mathematical statement is true. This question has a much more practical (and even cynical) counterpart: a student of a mathematics class wants to know what will the teacher accept as a correct solution of a homework problem. Here the logical structure of the claim is also very important, and we discuss several types of problems and their use in teaching mathematics.
\end{abstract}

%The organizers of the conference kindly invited me to speak about something related to constructivism and mathematical education. Indeed there is a theme that connects them, and I'd like to say a few words about it. First of all, a disclaimer: this is more a table talk than a research paper; I am grateful to the organizes for a chance to present some practical observations and general views.

\section{Language and truth}

Language is a mean of communication between people, somehow developed during the human history. It is natural to expect that its use should make human cooperation more efficient. Indeed, any complicated system (technical or natural) uses some kind of signals exchange between its parts. However, it is quite unclear from this viewpoint why indicative mood is used more often than the imperative mood. Why we often classify messages (sequences of signals) as being ``true'' or ``false'' and have even some special language tools (logical connectives and quantifiers) related to this classification? Moreover, sometimes even an order is presented as a statement (``You should do $X$.'' instead of ``Do $X$!'') Note that in programming the declarative style of programming and corresponding languages are much less popular.

\section{Mathematical statements and constructions}

Mathematicians and teachers of mathematics have to deal with a more specific problem: what is the meaning of a mathematical statement? Why and how can we classify some statement as true or false? The question becomes more difficult if the statement has a complicated logical structure (many alternating quantifiers) or deals with more abstract objects (like continuum hypothesis). 

The intuitionism/constructivism\footnote{We consider here only constructivism \emph{in mathematics}, not in any other sense (Wikipedia~\cite{wiki} lists $15$ different meanings, and only one is discussed here).} movement was an attempt to investigate different ways to understand mathematical statements on different levels of abstraction. One may pose the question in the following ``operational'' way: \emph{if we claim that some statement is true, what kind of mental construction should we have in mind to support this claim}? 

For example, how can we justify that $A\lor B$ is true? The natural answer: choose one of the statements $A$, $B$ and provide some justification for it. This natural answer, however, leads to non-classical logic where the law $A\lor\lnot A$ is not guaranteed (since we may not know whether $A$ is true or not).

\section{Understandable goals in teaching}

The question about operational definition of mathematical truth becomes quite practical from a student's perspective: \emph{given the problem, what should I show the teacher to get my work accepted}? Unfortunately, in many cases (I am afraid, in most of the cases, if we consider mathematics teaching in general) students just try to say random things in the hope they will satisfy the teacher, or try to perform some pattern matching using the examples shown during the previous class. This is quite frustrating (both for the students and the teacher) and, as a result, only a small fraction of intellectual abilities of students are used during a standard high school mathematics course (in most cases).

Long time ago, when Rubik's cube first appeared, many children learned how to assemble it --- not only  those that are good in mathematics, but also not so good ones. Obviously, the algorithm and its explanation is much more complicated (in terms of the number of cases, the necessity to explain 3D actions, the length of the description etc.) than \emph{any procedure or argument that appears in the high school mathematics program}. Why then is it more accessible? I think one of the reasons is that for the Rubik cube there is not need to guess what the teacher wants or be nice to her to get approved: \emph{you either are able to assemble the cube or not}, and that's all.

There are many similar stories: one attributed to Israel Gelfand says that a teacher for a long time tried to explain to his students how to compare fractions $2/3$ and $3/5$, failed and then --- as a last resort --- asked them what is better: three vodka bottles for five people, or two bottles for three. The answer was immediate: of course two for three is better.\footnote{An explanation: three for five is less than one bottle per person, so if someone come with his own bottle, the average increases and becomes four bottles for six people, i.e., two bottles for three.}

So why mathematics courses are so difficult? One of the reasons is that for students the mathematical classes look as a rather boring almost random trial-and-error guessing, and mathematical statements are just sequences of words and symbols, without any ``real'' meaning. No surprise that in this setting the learning is not very pleasant or effective. It seems that current (2023) AI learning systems exhibit a similar disappointing behavior when the topic is mathematics.

\section{Existential statements}

The lesson for mathematics teachers: one should choose problems that make sense for the students. After a short explanation the students need to be sure what the teacher wants and which kind of solutions she will accept. The easiest case is the case of \emph{existential} statements that claim that an object with some easily checkable properties exists. Then the solution consists of an example of such an object.

\begin{pr}
Find a  positive integer that becomes $57$ times smaller if we erase its first digit in its decimal representation.
\end{pr}

The solution $7125=57\cdot 125$ can be found by some algebraic reasoning, but this reasoning is not needed to get your solution accepted, an example is enough.

\begin{pr}
Draw a polygon and a point inside such that no side of the polygon is entirely visible from that point (each side is obscured partially by some other side).
\end{pr}

One of possible solutions is shown at the picture: long sides are partially obscured while the short sides are not visible at all.
\begin{center}
\includegraphics[width=0.3\textwidth]{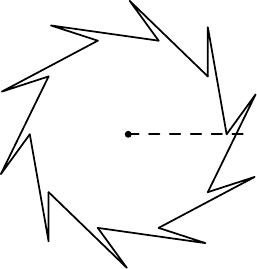}
\end{center}
Again it is pretty obvious both for students and for teachers whether a given polygon solves the problem. (There could be minor troubles with ``self-intersecting polygons''.)

\begin{pr}\label{positive-negative}
Write $10$ numbers in a line in such a way that sum of every three neighbors is positive, and the sum of all ten numbers is negative.
\end{pr}

One of the possible answers:
\[
{-}9\quad  5\quad  5\quad  {-9}\quad  5\quad   5\quad  {-9}\quad   5\quad   5\quad   {-9}.
\]
Here every three neighbors have sum $1$, and all the ten numbers have sum $-6$.
\medskip

The problem may be not purely mathematical:
\begin{pr}
Given an A4 sheet of paper and scissors, cut a hole in the paper that is large enough so you can go through it.
\end{pr}
Even if this problem is not purely mathematical, it is quite clear what kind of object you should provide to solve it:  something with a hole big enough so that you can put it around yourself. (There is a myth where Didone's solution of a similar problem turned out to be convincing.)
\begin{center}
\includegraphics[width=0.6\textwidth]{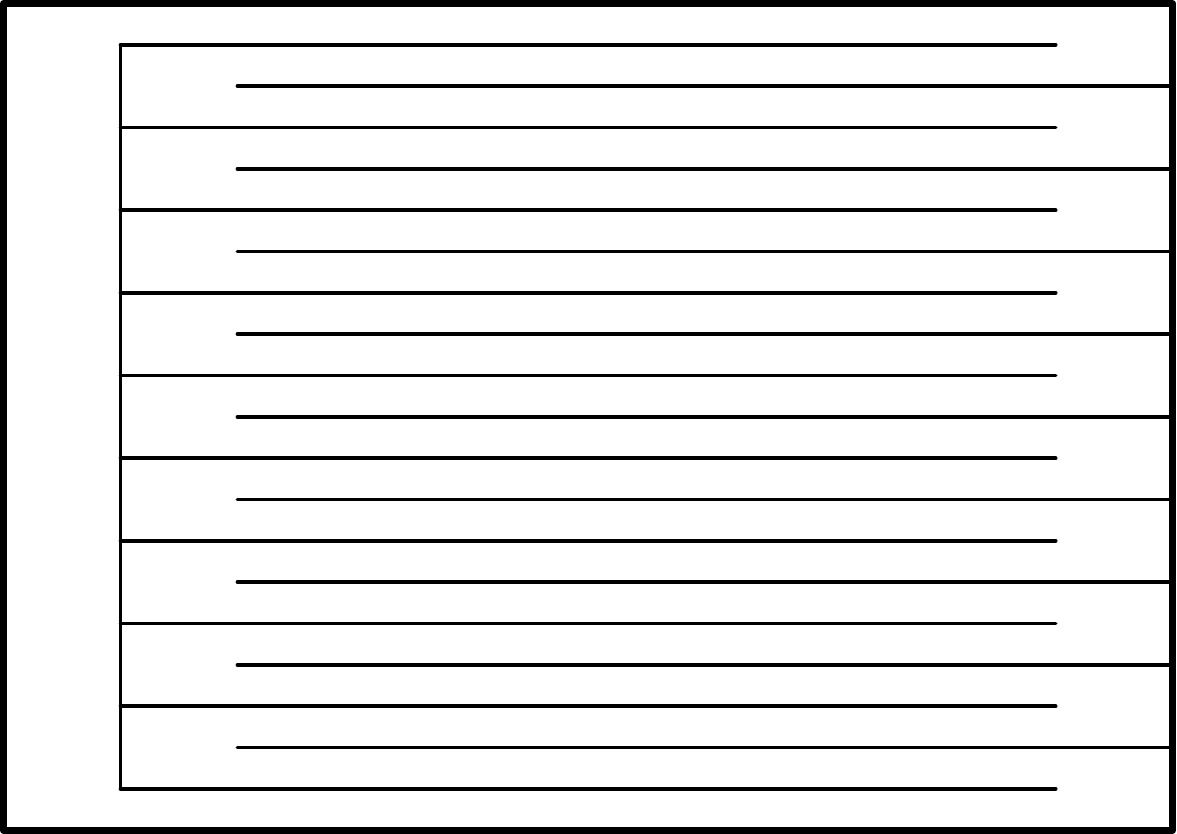}
\end{center}

Let me stress again the main point: the existential nature of all these problems makes them suitable for teaching: the construction presented as a solution can be easily checked without a teacher. (Would the Sudoku game be popular if you had to present the solution to someone to find out whether it is valid? Probably not.) 

They are also good for mathematical competitions: for grading it is enough to look at the answer, there is no need to follow long (and usualy badly written) arguments. Of course, one can choose any problem with a numerical answer and just ask to provide the answer without explanations --- but then participant will not understand why her answer is rejected, and this does not provide the feeling of fair grading.

\section{Universal statements}

Universal statements (at least from a classical viewpoint) are negations of existential statements. Let us consider a problem quite similar to problem~\ref{positive-negative}.

\begin{pr}
Put $10$ numbers \emph{along a circle} in such a way that the sum of every three neighbors is positive and the sum of all ten numbers is negative.
\end{pr}

The difference is that now the numbers are placed around the circle (this creates new sums that should also be positive). This change makes the task impossible. Indeed, there are $10$ sums (of three numbers each) that should be positive, and if we sum up all these $10$ sums, we get a sum where each number appears three times, so it is the sum of all numbers multiplied by three, and at the same time should be positive and negative --- a contradiction.

So the task is unsolvable, and this unsolvability is a universal statement. The operational meaning of the universal statement is less clear: what kind of argument should we provide to justify it? It is a difficult question (both for constructivists and math teachers). Still it is possible to make this question more practical. If somebody claims that the task is impossible, we can suggest her to make a bet (according to the proverb that says ``put your money where your mouth is''). Often this question changes the attitude of a student completely: before it was just an attempt to guess what answer the teacher expects, and suddenly it becomes a very practical question (especially for people who have lost similar bets earlier).

Still what kind of arguments should the student provide to justify her answer? As na\"\i ve constructivists, we say that any arguments are OK if they are convincing. Imagine that Alice wants to convince Bob that he can safely make a bet (with Charlie). Alice may use whatever arguments she wants --- but they should really convince Bob.  

\begin{pr}
Cut an $8\times 8$ board without two opposite corners into $1\times 2$ domino tiles.
\end{pr}
\begin{center}
\includegraphics[width=0.3\textwidth]{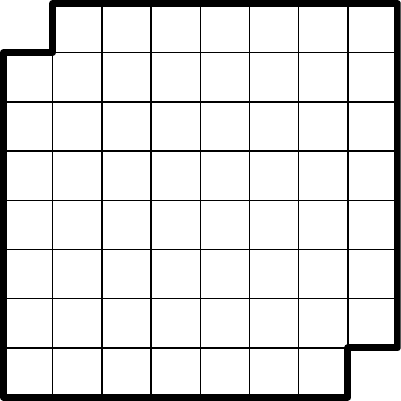}
\end{center}

Most students easily find that their attempts to do this are unsuccessful. But is this negative experience enough for them to safely make a bet? May be the malicious opponent had some ingenious way of cutting the board? Is there some reason to be sure that this is not possible? Indeed, we may use the chessboard coloring 
\begin{center}
\includegraphics[width=0.3\textwidth]{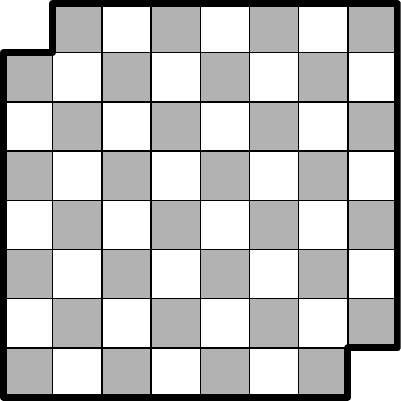}
\end{center}
with $32$ black cells and $30$ white cells and note the each domino consists of one black and one white cell. If a cut were possible, then the numbers of white and black cells would be equal (being equal to the number of dominoes) --- and this is not the case.

\section{Maximization}

There is one more class of problem that naturally combines existential and universal statements.
\begin{pr}
Put a maximal number of knights on the chessboard in such a way that none of them attack any other.
\end{pr}
The quantifier prefix here is $\exists\forall$ (there exist some placement of knights such that no other placement has more knights). The goal of the student is therefore twofold:
\begin{itemize}
\item show a position with some number of knights not attacking each other;
\item convince someone that she may safely make a bet saying that no one would come with a position with more knights not attacking each other.
\end{itemize}
The solution is easy to explain (though, maybe, not so easy to invent). Here is a position with $32$ knights
\fenboard{N1N1N1N1/1N1N1N1N/N1N1N1N1/1N1N1N1N/N1N1N1N1/1N1N1N1N/N1N1N1N1/1N1N1N1N w - - 0 20}
\begin{center}
\showboard
\end{center}
Note that any two knights on white cells cannot attack each other.

How can we be sure that no other valid position contains more than $32$ knights? Let us split all the cells on a chessboard into pairs. Namely, we first cut the board into $8$ rectangles of size $4\times 2$, and then group $8$ cells in each rectangle into four pairs:
\begin{center}
\includegraphics[width=0.2\textwidth]{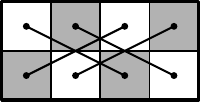}
\end{center}
Each pair may contain at most one knight (otherwise two knights attack each other), therefore the total number of knights in a valid configuration is bounded by the number of pairs, i.e.,~$32$.

\section{More complicated statements}

The classical mathematicians (and teachers of mathematics) then try to create a \emph{psychologically robust illusion} that more complicated mathematical statements also have some meaning in the sense that they are somehow objectively true or false. (Of course, on the next level mathematical logicians will explain why it is only an illusion, but as a first approximation it works.) There are some tricks that help the teachers to achieve this goal.

\subsection{Reinterpretation}

Alternating quantifiers can be interpreted in game terms. For example, what construction a student should provide to show that $\lim_{n\to\infty} a_n =0$ for some sequence $a_n$? The definition says
\[
(\forall \varepsilon>0)\, \exists N\, (\forall n>N) [ |a_n| < \varepsilon]. 
\]
So, says the teacher, imagine that we play a game: I make the first move and choose some positive $\varepsilon$, then you may reply by choosing some $N$, then I choose some $n>N$, and then we look together who is the winner: you win if $|a_n|$ is smaller than $\varepsilon$. Are you sure that you can with this game? (If the student says yes, one can then try to beat her in this game and see whether this belief is based on some specific strategy.)

In logical terms, we reinterpret a statement with long quantifies prefix as a statement about existence of a strategy that wins the game --- and this is a statement of the $\exists\forall$-type. Of course, the word ``strategy'' mus be understood more or less constructively (as an ``arbitrary function'', an algorithm, etc.)

N.N.~Konstantinov noted that even in the golden years of \emph{mekhmat} (Moscow State Lomonosov University mathematics department) quite a few students were unable to win such a game (for a similar statement: the continuity of a square root function at some point) even after passing calculus exam. (A bad symptom for teaching quality.)

\subsection{Intermediate notions}

Another way to become accustomed to long quantifier prefixes is to split them using intermediate definitions. Let me explain what I mean using the same toy example (the definition of limit). This topic is traditionally hard, especially for high school students, so the following approach is used in the nice book of Alexander Kirillov~\cite{kirillov}. 

Let us call an interval  $I$ a \emph{trap} for a sequence $a_0,a_1,\ldots$ if all $a_n$ with sufficiently large $n$ belong to $I$. (Imagine that $a_n$ is a position of an object at moment $n$.) Then one could start solving problems that use this notion. This is easier, since we have one quantifier less than in the definition of a limit. For example, one may ask why an intersection of two traps for the same sequence is also a trap for it, or ask to find a trap for $10a_n$ if $(1,2)$ is a trap for $a_n$, etc.

After this notion becomes familiar, the requirement in the definition of a limit sounds quite simple, only one universal quantifier remains: \emph{$(-\varepsilon,\varepsilon)$ is a trap for every~$\varepsilon$}.

\subsection{Possible constructions as finished ones}

\begin{pr}
Prove that some number of the form $111\ldots 1$ is divisible by $179$.
\end{pr}

The following solution looks convincing, though it does not provide such a number explicitly. Consider the numbers $1, 11, 111, 1111, \ldots$ (up to $1$ repeated $179$ times) modulo $179$. If none of them is a multiple of $179$, then some two should have the same remainder modulo $179$, so their difference ($11\ldots 1100\ldots00$) is divisible by $179$, and remains divisible after the trailing zeros are deleted.

Why this argument looks convincing? One may say that though we do not show such a number, we describe a clear (though boring) procedure to find it. This is quite standard practice: when asked how many $32$-bit strings exists, we answer ``$2^{32}$''. But, strictly speaking, we do not provide the answer in this way, we just explain what should be done to find it (we need to multiply $32$ numbers each equal to $2$). Even a more explicit answer $4\,294\,967\,296$ also just encodes some procedure: take $4$, multiply it by $10$ nine times, then add $2$ multiplied by $10$ eight times, etc.

One could say that we agree to think about some described processes as finished ones. Of course, in this case this leap does not look dangerous, but then it is easy to apply the same intuition to infinite processes (``look at all numbers greater than $N$; if they are no twin primes among them, then do something, if there are, do something else'', etc.)

\subsection{Is it needed?}

All these examples lead to some natural question: why do we try to create an illusion for our students if we know that it is only an illusion? Wouldn't it be better to start teaching honest constructive mathematics without all these dirty tricks? Do we need to learn the classical results before trying the extract constructive statements from its proofs or provide new constructive proofs? (A good example of the latter is the Bishop--Vyugin constructive proof of the ergodic theorem~\cite{vyugin}.) 

The experience of mathematical community seems to suggest a positive answer: for some reasons it is indeed easier to start with ideal illusory objects before trying to translate the arguments to the ground level, and teach the ``classical'' mathematics before the ``constructive'' one.
\bigskip

But why?

\end{document}